\documentclass{amsart}

\usepackage{amssymb}
\usepackage[all]{xypic}
\xyoption{dvips}

\theoremstyle{definition}
\newtheorem{ntn}{Notation}

\theoremstyle{plain}
\newtheorem{lem}[ntn]{Lemma}

\newtheorem{thm}[ntn]{Theorem}

\theoremstyle{remark}

\newcommand{\m}{\mathfrak{m}}

\newcommand{\C}{\mathbb{C}}
\newcommand{\N}{\mathbb{N}}
\renewcommand{\O}{\mathcal{O}}
\newcommand{\Q}{\mathbb{Q}}
\newcommand{\wh}{\widehat}
\newcommand{\xymat}{\SelectTips{cm}{}\xymatrix}

\DeclareMathOperator{\Ann}{Ann}
\DeclareMathOperator{\Der}{Der}
\DeclareMathOperator{\Ext}{Ext}
\DeclareMathOperator{\Hom}{Hom}
\DeclareMathOperator{\res}{res}
\DeclareMathOperator{\Sing}{Sing}
\DeclareMathOperator{\tr}{tr}

\begin{document}

\title[Quasihomogeneity of isolated hypersurface singularities]{
Quasihomogeneity of isolated hypersurface singularities and logarithmic cohomology}
\author{Michel Granger and Mathias Schulze}

\begin{abstract}
We characterize quasihomogeneity of isolated hypersurface singularities by the injectivity of the map induced by the first differential of the logarithmic differential complex in the top local cohomology supported in the singular point.
\end{abstract}

\thanks{The second author gratefully acknowledges financial support from the Humboldt foundation.}



\maketitle

\section{Formulation of the result}

We consider the germ of an isolated hypersurface singularity $D\subseteq(\C^n,0)$.
Let $\O=\O_{\C^n,0}$ be the ring of germs of holomorphic functions on $\C^n$ at the origin with maximal ideal $\m$ and $f\in\m^2$ a reduced defining equation of $D$.
For notational simplicity, we shall write $\C^n$ instead of $(\C^n,0)$.
In terms of coordinates $x_1,\dots,x_n$ on $\C^n$, $\O=\C\{x_1,\dots,x_n\}$ is the ring of convergent power series and $\m=\langle x_1,\dots,x_n\rangle$.
We abbreviate by $f_i=\frac{\partial f}{\partial{x_i}}$, $i=1,\dots,n$, the partial derivatives of $f$.

Let $\Omega^\bullet=\Omega_{\C^n,0}^\bullet$ be the complex of holomorphic differential forms at the origin and denote by $\Omega^\bullet(\star D)=\Omega^\bullet[f^{-1}]$ the complex of meromorphic differential forms with poles along $D$.
K.~Saito \cite[\S 1]{Sai80} introduced the subcomplex of logarithmic differential forms $\Omega^\bullet(\log D)\subseteq\Omega^\bullet(\star D)$ along $D$.
It is the maximal subcomplex consisting of forms with at most a simple pole along $D$.
Let $\Der=\Der_\C\O$ be the module of $\C$-linear derivations of $\O$.
It is a free $\O$-module of rank $n$ with a basis formed by the partial derivatives $\partial_1,\dots,\partial_n$ with respect to $x_1,\dots,x_n$.
The module of logarithmic vector fields $\Der(-\log D))\subseteq\Der$ consists of all vector fields tangent to the smooth part of $D$.
Under the natural inner product, $\Der(-\log D))$ and $\Omega^1(\log D)$ are mutually $\O$-dual.

Our object of interest in this article is the natural map
\[
\xymat{
H^n_0(\O)\ar[r]^-{d_1}&H^n_0(\Omega^1(\log D))
}
\]
induced by the first differential of $\Omega^\bullet(\log D)$ in the top local cohomology supported at the origin. 
We shall show that its injectivity is equivalent to the quasihomogeneity of the singularity.

\begin{thm}\label{0}
$D$ is quasihomogeneous if and only if $d_1$ is injective.
\end{thm}

For non--isolated singularities, at least one implication fails:
The divisor $D\subseteq\C^3$ defined by $f=xy(x+y)(x+yz)$ is not quasihomogeneous but $d_1$ is injective.

\section{Relation to the logarithmic comparison theorem}

The map $d_1$ is well-known in the context of the logarithmic comparison theorem:
By Grothendieck's comparison theorem \cite{Gro66}, the cohomology of the complement $U=\C^n\smallsetminus D$ of a germ of a reduced divisor $D\subseteq\C^n$ is computed by the complex $\Omega^\bullet(\star D)$.
If the inclusion $\Omega^\bullet(\log D)\subseteq\Omega^\bullet(\star D)$ is a quasi--isomorphism then $\Omega^\bullet(\log D)$ computes the cohomology of $U$ and one says that the logarithmic comparison theorem, or briefly LCT, holds for $D$.
The characterization of LCT is an open problem and the subject of active current research \cite{Tor05}.
Mainly the two extremal cases of isolated singularities and free divisors have been studied.
By definition, $D$ is free if $\Der(-\log D)$ is $\O$-free of rank necessarily equal to $n$.
Freeness of $D$ means that it is either smooth or its singular locus $\Sing(D)$ is of pure codimension $1$ in $D$.
Thus the two cases coincide for singularities in dimension $n=2$ and are disjoint if $n\ge3$.

By \cite[Sect.~2]{CNM96,CMNC02}, the map $d_1$ occurs as the differential $d_1^{0,n-1}$ in the $1$-page 
\[
E_1^{p,q}=H^q(\C^n\smallsetminus\{0\},\Omega^p(\log D))
\]
of a spectral sequence $E$ which converges to the cohomology of $U$ under the assumption that LCT holds for $D\smallsetminus\{0\}\subseteq\C^n\smallsetminus\{0\}$.
In particular, this assumption is fulfilled in the case of isolated singularities.
Since sections of $\Omega^\bullet(\log D)$ defined outside $\Sing(D)$ extend uniquely to $\C^n$, the $q=0$ row of $E_1$ equals $\Omega^\bullet(\log D)$.
Thus, if LCT holds also at the origin, the spectral sequence $\bar E$ with $1$-page $\bar E_1$ obtained from $E_1$ by removing the $q=0$ row converges to zero.

In the case of a free divisor $D$, $\Omega^\bullet(\log D)$ consists of free $\O$-modules and hence $E_1^{p,q}=0$ for $q\ne0,n-1$.
This implies that the injectivity of $d_1$ is a necessary condition for LCT  \cite[Sect.~3]{CMNC02}.
In dimension $n=2$, injectivity of $d_1$ is even equivalent to LCT and to quasihomogeneity \cite[Thm.~3.5]{CMNC02}.

If $f\in\m_xJ_f$ at any point $x\in D$ for a reduced defining equation $f\in\m$ of $D$ then $D$ is called strongly Euler homogeneous.
In the case of free divisors, the main conjecture on LCT states that LCT implies strong Euler homogeneity \cite[Conj.~1.4]{CMNC02}.
We have proved this conjecture in dimension $n=3$ \cite{GS06}.

For isolated singularities, strong Euler homogeneity reduces to quasihomogeneity by \cite{Sai71}.
M.~Holland and D.~Mond \cite{HM98} have given an explicit characterization of LCT for quasihomogeneous isolated singularities $D$:
Let $f\in\m^2$ be a reduced defining equation of $D$ which is quasihomogeneous with respect to the necessarily strictly positive weights $w_1,\dots,w_n\in\Q$.
Then LCT holds for $D$ if and only if 
\[
(\O/\langle f_1,\dots,f_n\rangle)_{k-\sum_{i=1}^nw_i}=0\text{ for }k=1,\dots,n-2.
\]
In particular, quasihomogeneity does not imply LCT for isolated singularities.
On the other hand, we do not know a counter--example to the opposite implication.
Our result possibly reduces the complete characterization of LCT for isolated singularities to the statement:
LCT implies injectivity of $d_1$ for isolated singularities.
But it is not clear if the injectivity of $d_1$ contradicts to $\bar E$ converging to zero.

\section{Proof of the theorem}

We denote the $n$th local cohomology module of $\O$ supported at the origin by
\[
H:=H_0^n(\O)
=H^{n-1}(\C^n\smallsetminus\{0\},\O)
=\frac{\O[\frac{1}{x_1\cdots x_n}]}{\sum_{i=1}^n\O[\frac{1}{x_1\cdots\wh{x_i}\cdots x_n}]}.
\]
By right exactness of the tensor product $\otimes:=\otimes_\O$, this implies
\[
H_0^n(\Omega^1(\log D))=H^{n-1}(\C^n\smallsetminus\{0\},\Omega^1(\log D))=\Omega^1(\log D)\otimes H.
\]
The total differential $d:\O\to\Omega^1(\log D)$ induces a natural map
\[
d_1:H=H_0^n(\O)\to H_0^n(\Omega^1(\log D))=\Omega^1(\log D)\otimes H
\]
in local cohomology.
The perfect pairing $\Der(-\log D)\otimes\Omega^1(\log D)\to\O$ \cite[Lem.~1.6]{Sai80} identifies $\Omega^1(\log D)$ to $\Hom(\Der(-\log D),\O)$.
Applying $\otimes H$, it induces a pairing
\begin{equation}\label{2}
\langle,\rangle:\Der(-\log D)\otimes H_0^n(\Omega^1(\log D))=\Der(-\log D)\otimes(\Omega^1(\log D)\otimes H)\to H.
\end{equation}
such that $\langle\delta,(\omega\otimes[g])\rangle=\langle\delta,\omega\rangle [g]$.
In particular, $d_1([g])=\sum_{i=1}^ndx_i\otimes[\partial_i(g)]$ for $g\in\O[\frac{1}{x_1\cdots x_n}]$ and hence
\begin{equation}\label{3}
\langle\delta,d_1([g])\rangle=\sum_{i=1}^n\langle\delta,dx_i\rangle[\partial_i(g)]=\delta([g]).
\end{equation}  
Denoting $\Hom=\Hom_\O$, the pairing (\ref{2}) corresponds to the natural map 
\begin{equation}\label{5}
\xymat{
\Hom(\Der(-\log D),\O)\otimes H=\Omega^1(\log D)\otimes H\ar[r]^-\pi&\Hom(\Der(-\log D),H).
}
\end{equation}
According to \cite[Lem.~1.11]{HM98}, $\Omega^1(\log D)=\O\cdot\frac{df}{f}+\Omega ^1$ and there is an exact sequence
\begin{equation}\label{4}
\xymat{
0\ar[r]&\O\ar[r]^-{\alpha^t}&\O^{n+1}\cong\O\oplus\Omega^1\ar[r]^-\beta&\Omega^1(\log D)\ar[r]&0
}
\end{equation}  
where $\alpha=(f,f_1,\dots,f_n)$ and $\beta=(\frac{df}{f},-dx_1,\dots,-dx_n)$.
Dualizing against $\O$ yields
\[
\xymat{
0\ar[r]&\Der(-\log D)\ar[r]^-{\beta^t}&\O^{n+1}\ar[r]^-\alpha&J\ar[r]&0,\quad J:=\langle f,f_1,\dots,f_n\rangle,
}
\]
and $\Ext^1(\Omega^1(\log D),\O)\cong\O/J$.
Dualizing against $\O$ and $H$ yields the sequences
\begin{gather}
\label{6}
\xymat@C=22pt{
0\ar[r]&\Hom(J,\O)\ar[r]&\O^{n+1}\ar[r]^-{(\beta^t)^t}&\Hom(\Der(-\log D),\O)\ar[r]&0,
}\\
\label{7}
\xymat@C=14pt{
0\ar[r]&\Hom(J,H)\ar[r]&H^{n+1}\ar[r]&\Hom(\Der(-\log D),H)\ar[r]&\Ext^1(J,H)\ar[r]&0.
}
\end{gather}
By reflexivity of $\Der(-\log D)$ \cite[Cor.~1.7]{Sai80}, the sequences (\ref{4}) and (\ref{6}) coincide.
In particular, $(\beta^t)^t=\beta$ is surjective and $\Hom(J,\O)=\O$ as the kernel of $\beta$ which is directly obvious.
From sequences (\ref{6}) and (\ref{7}), we obtain the commutative diagram 
\begin{equation}\label{8}
\xymat{
\Hom(J,\O)\otimes H \ar[r]\ar[d]&H^{n+1} \ar@{->>}[r]^-{\beta\otimes1}\ar@{=}[d]&\Hom(\Der(-\log D),\O)\otimes H\ar[d]^-\pi\\
\Hom(J,H)\ar@{^(->}[r]^-\iota&H^{n+1}\ar[r]^-\gamma&\Hom(\Der(-\log D),H)
}
\end{equation}
involving the map $\pi$ in (\ref{5}).

\begin{lem}\label{1}
If $D$ is not Euler homogeneous then $\tr\delta_0=0$ for all $\delta\in\Der(-\log D)$.
In particular, $\pi(d_1([\frac{1}{x_1\dots x_n}]))(\delta)=\langle\delta,d_1([\frac{1}{x_1\dots x_n}])\rangle=0$ for all $\delta\in\Der(-\log D)$.
\end{lem}

\begin{proof}
By the formal structure theorem \cite[Thm.~5.4]{GS06}, there is a system of generators $\mathcal S$ of the $\m$-adic completion $\wh\Der(-\log D)$ of $\Der(-\log D)$ such that, for $\delta\in\mathcal S$, either $\delta_0$ is nilpotent and hence $\tr\delta_0=0$ or $\delta=\delta_0$ is semisimple and $\delta(\hat f)\in\Q\cdot\hat f$ where $\hat f:=\hat u\cdot f$ for some unit $\hat u\in\wh\O^\star$.
In the second case, $\delta\in\Ann_{\wh\Der}(\hat f)$ if $D$ and hence $\hat f$ is not Euler homogeneous.
Since $D$ is an isolated singularity,
\[
\Ann_{\wh\Der}(\hat f)=\langle\hat f_i\partial_j-\hat f_j\partial_i\mid 1\le i<j\le n\rangle,\quad\hat f_i:=\partial_i(\hat f).
\]
By a formal coordinate change \cite[Lem.~4.1]{Arn73}, $\hat f=x_1^2+\cdots+x_k^2+g(x_{k+1},\dots,x_n)$ with $g\in\wh\m^3$ and $\tr\delta_0=0$ follows also in the second case from
\[
\delta_0\in\Ann_{\wh\Der}(\hat f)_0=\langle x_i\partial_j-x_j\partial_i\mid 1\le i<j\le k\rangle.
\]
Since $\Der(-\log D)_0=\wh\Der(-\log D)_0$ and by additivity and coordinate invariance of $\tr\delta_0$, this proves the first claim.
The second follows from (\ref{3}) and the obvious formula $\delta([\frac{1}{x_1\dots x_n}])=\tr(\delta_0)[\frac{1}{x_1\dots x_n}]$ for any $\delta\in\m\cdot\Der$.
\end{proof}

We continue our study of diagram (\ref{8}) under the assumption that $D$ is not Euler homogeneous.
By Lemma \ref{1}, $\pi(d_1([\frac{1}{x_1\cdots x_n}]))=0$ and $d_1([\frac{1}{x_1\cdots x_n}])=(\beta\otimes1)(\psi)$ where $\psi=(0,u_1,\dots,u_n)\in H^{n+1}$ and $u_i=-\partial_i([\frac{1}{x_1\cdots x_n}])=[\frac{1}{x_1\cdots x_i^2\cdots x_n}]$.
Then $\gamma(\psi)=0$ and hence $\psi=\iota(\varphi)$ where $\varphi\in\Hom(J,H)$ is characterised by
\begin{equation}\label{9}
\varphi(f_i)=\left[\frac{1}{x_1\cdots x_i^2\cdots x_n}\right],\quad\varphi(f)=0.
\end{equation}
By injectivity of $\iota$, $d_1([\frac{1}{x_1\cdots x_n}])=0$ if and only if $\varphi$ extends to $\O$.

\begin{lem}\label{10}
The map $\varphi\in\Hom(J,H)$ which fulfils (\ref{9}) extends to a $\tilde\varphi\in\Hom(\O,H)$.
\end{lem}

\begin{proof}
Denote by $A:=\C[x_1,\dots,x_n]\subseteq\O$ the polynomial ring and by $M^\star:=\Hom_\C(M,\C)$ the algebraic dual of an $A$- or $\O$-module $M$.
The residue pairing \cite[Sect.~B.1]{Ems78} $\langle,\rangle:H\otimes_\C\O\to\C$ is defined by $\langle[g],h\rangle:=\res([g\cdot h])$ where the residue map $\res:H\to\C$ sends $[g]\in H$ to the coefficient of $\frac{1}{x_1\cdots x_n}$ in $g$.
More explicitly,
\[
\left\langle\left[\sum_{i\in\N^n}\frac{a_i}{x_1^{i_1+1}\cdots x_n^{i_n+1}}\right],\sum _{i\in\N^n}b_ix_1^{i_1}\cdots x_n^{i_n}\right\rangle=\sum_{i\in\N^n}a_ib_i.
\]
The residue pairing defines an injection of $H$ into $\O^\star$.
Restricted to $A$, the pairing injects $H$ into $A^\star=\C[[x_1^{-1},\dots,x_n^{-1}]]\cdot\frac{1}{x_1\cdots x_n}$.
By (\ref{9}), $\varphi$ annihilates the ideal $K:=\langle f\rangle+\m^2\cdot\langle f_1,\dots,f_n\rangle$ and induces a $\bar\varphi\in\Hom(J/K,H)$.
Since $K$ has finite colength, $\m^{N+1}\subseteq K$ for some $N>0$ and hence the finite vector space
\[
K^\perp\subseteq(\m^{N+1})^\perp=\bigoplus_{|i|\le N}\frac{\C}{x_1^{i_1+1}\cdots x_n^{i_n+1}}
\]
identifies to $(K\cap A)^{\perp}\subseteq A^\star$.
Then the residue pairing induces an isomorphism
\[
\ell:K^{\perp}=(K\cap A)^{\perp}\to(A/K\cap A)^{\star}=(\O/K)^{\star}
\] 
defined by $\ell([g])([h])=\langle[g],h\rangle=\res([g\cdot h])$ for $[g]\in K^{\perp}$ and $[h]\in\O/K$.
Thus, the linear form $\res\circ\bar\varphi\in(J/K)^\star$ extends to $\ell([g])$ for some $[g]\in K^{\perp}$ which means that $\res(\varphi(h))=\res([g\cdot h])$ for all $h\in J$.
Now, define $\tilde\varphi\in\Hom(\O,H)$ by $\tilde\varphi(h):=[g\cdot h]$ for all $h\in\O$.
Then, for $h\in J$,
\[
\langle\tilde\varphi(h),x^\alpha\rangle=\langle[g\cdot h],x^\alpha\rangle=\res([g\cdot x^\alpha h])=\res(\varphi(x^\alpha h))=\langle\varphi(h),x^\alpha\rangle
\]
for all $\alpha\in\N^n$ and hence $\tilde\varphi(h)=\varphi(h)$ which implies that $\tilde\varphi$ is an extension of $\varphi$.
\end{proof}

\begin{proof}[Proof of Theorem \ref{0}]
For isolated singularities, Euler and quasihomogeneity are equivalent by \cite{Sai71} and $h^p(H^n_0(\Omega^\bullet(\log D)))=0$ for all $p$ by \cite[Cor.~1.8]{HM98}.
For non Euler homogeneous $D$, $0\ne[\frac{1}{x_1\cdots x_n}]\in\ker d_1$ by Lemma \ref{10}. 
\end{proof}

\bibliographystyle{amsalpha}
\bibliography{isd1}

\end{document}